\documentclass[a4paper]{amsart}
\usepackage{multicol,ifthen}
\usepackage{amsmath}
\usepackage{amssymb}
\usepackage{amsthm} 
\usepackage{enumerate}

\begin{document}

\newcommand{\F}{\mathcal{F}}
\newcommand{\R}{\mathbb R}
\newcommand{\T}{\mathbb T}
\newcommand{\N}{\mathbb N}
\newcommand{\Z}{\mathbb Z}
\newcommand{\C}{\mathbb C}  
\newcommand{\h}[2]{\mbox{$ \widehat{H}^{#1}_{#2}(\R)$}}
\newcommand{\hh}[3]{\mbox{$ \widehat{H}^{#1}_{#2, #3}$}} 
\newcommand{\n}[2]{\mbox{$ \| #1\| _{ #2} $}} 
\newcommand{\x}{\mbox{$X^r_{s,b}$}} 
\newcommand{\xx}{\mbox{$X_{s,b}$}}
\newcommand{\X}[3]{\mbox{$X^{#1}_{#2,#3}$}} 
\newcommand{\XX}[2]{\mbox{$X_{#1,#2}$}}
\newcommand{\q}[2]{\mbox{$ {\| #1 \|}^2_{#2} $}}
\newcommand{\e}{\varepsilon}
\newcommand{\om}{\omega}
\newcommand{\lb}{\langle}
\newcommand{\rb}{\rangle}
\newcommand{\ls}{\lesssim}
\newcommand{\gs}{\gtrsim}
\newcommand{\pd}{\partial}
\newtheorem{lemma}{Lemma} 
\newtheorem{kor}{Corollary} 
\newtheorem{theorem}{Theorem}
\newtheorem{prop}{Proposition}

\title[Wave equation with quadratic nonlinearities in 3 D]{On the wave equation with quadratic nonlinearities in three space dimensions}

\author[Axel~Gr{\"u}nrock]{Axel~Gr{\"u}nrock}

\address{Axel~Gr{\"u}nrock: Rheinische Friedrich-Wilhelms-Universit\"at Bonn,
Mathematisches Institut, Endenicher Allee 60, 53115 Bonn, Germany.}
\email{gruenroc@math.uni-bonn.de}

\thanks{The author was partially supported by the Deutsche Forschungsgemeinschaft, Sonderforschungsbereich 611.}

\subjclass[2000]{35L70}

\begin{abstract}
 The Cauchy problem for the nonlinear wave equation
$$\Box u=(\partial u)^2, \qquad u(0)=u_0, \,\,u_t(0)=u_1$$
in three space dimensions is considered. The data $(u_0,u_1)$ are assumed to belong to
$\widehat{H}^r_s(\R^3) \times \widehat{H}^r_{s-1}(\R^3)$, where $\widehat{H}^r_s$ is defined
by the norm
$$\n{f}{\widehat{H}^r_s} := \n{\langle \xi \rangle ^s\widehat{f}}{L^{r'}_{\xi}},\quad \langle \xi \rangle=(1+|\xi|^2)^{\frac12}, \quad \frac{1}{r}+\frac{1}{r'}=1.$$
Local well-posedness is shown in the parameter range $2 \ge r >1$, $s > 1 + \frac{2}{r}$. For $r=2$ this coincides with the
result of Ponce and Sideris, which is optimal on the $H^s$-scale by Lindblad's counterexamples, but nonetheless leaves a gap
of $\frac12$ derivative to the scaling prediction. This gap is closed here except for the endpoint case. Corresponding results
for $\Box u = \partial u^2$ are obtained, too.
\end{abstract}

\keywords{Nonlinear wave equation -- Cauchy Problem -- local wellposedness -- generalized Bourgain spaces}

\maketitle

\section{Introduction and main results}

In this note we consider the Cauchy problem
\begin{equation}
 \label{CP}
\Box u=(\partial_t^2 - \Delta)u= B_k(u,u), \qquad u(0)=u_0, \,\,u_t(0)=u_1
\end{equation}
for the nonlinear wave equation in $\R^3$, where the right hand side is given by
\begin{equation}
 \label{B12}
B_1(u,v)=\partial(uv) \qquad {\mbox{or}} \qquad B_2(u,v)=\partial u \partial v
\end{equation}
with $\partial \in \{\partial_t , \partial_{x_1}, \partial_{x_2}, \partial_{x_3}\}$, and no special structure of the bilinear forms $B_k$, $k \in \{1,2\}$, such as a null-structure is assumed. Concerning the local well-posedness (LWP) of this problem with data
$(u_0,u_1) \in H^s(\R^3)\times H^{s-1}(\R^3)$ the following is known. For $s>k + \frac12$ energy estimates can be applied to obtain
an affirmative result. Ponce and Sideris showed in \cite{PS93} how to improve this down to $s>k$ by using Strichartz inequalities.
Further progress is possible, if the nonlinearity satisfies a null-condition such as
$$\widetilde{B}_2(u,v)=\lb\nabla_xu, \nabla_x v\rb- \partial_t u \partial_t v,$$
see the work of Klainerman and Machedon \cite{KM95}, \cite{KM96}, \cite{KM97}, who used wave Sobolev spaces to exploit the
null-structure of the bilinear terms, thus reaching LWP for $s>s_c(k)=k-\frac12$, which is here the critical Sobolev regularity by
scaling considerations. If no such structure is present in the quadratic term, one has in fact ill-posedness of the Cauchy
problem \eqref{CP} for $s \le k$, as the sharp counterexamples of Lindblad show, see \cite{Lb93}, \cite{Lb96}, \cite{Lb98}. So in
general there is a gap of half a derivative between the optimal LWP result on the $H^s$-scale and the scaling prediction.

\qquad

For several important nonlinear dispersive equations in \emph{one} space dimension - such as cubic NLS and DNLS, KdV, mKdV and its
higher order generalizations - there is a similar gap between the best possible LWP result in $H^s$ and the critical regularity.
In the case of cubic nonlinearities this can be closed almost completely by considering data in the spaces $\widehat{H}^r_s$,
defined by the norm
$$\n{f}{\widehat{H}^r_s} := \n{\langle \xi \rangle ^s\widehat{f}}{L^{r'}_{\xi}},\quad \langle \xi \rangle=(1+|\xi|^2)^{\frac12}, \quad \frac{1}{r}+\frac{1}{r'}=1,$$
see \cite{G04}, \cite{G05}, \cite{GV}, \cite{G09}; for an application in the periodic setting cf. \cite{GH}. The purpose of this note is to show that the methods developed in the one-dimensional framework also apply to the nonlinear wave equation \eqref{CP},
\eqref{B12} in three space dimensions and give LWP for data $(u_0,u_1) \in \widehat{H}^r_s(\R^3) \times \widehat{H}^r_{s-1}(\R^3)$,
provided $1<r\le 2$ and $s>s_k(r):=k+1-\frac{2}{r'}$. In the limit $r \rightarrow 1$ we can almost reach the space
$\widehat{H}^1_{k+1}(\R^3) \times \widehat{H}^1_{k}(\R^3)$, which is critical by scaling. To prove this result we use an appropriate variant of Bourgain's Fourier restriction norm method, see \cite[Section 2]{G04}, and estimates for products of two free solutions
of half-wave equations. The latter are very much in the spirit of the work of Foschi and Klainerman \cite{FK} and can be seen as
bilinear substitutes and refinements of the Strichartz inequalities for the three-dimensional wave equation.

\section{General arguments, function spaces, and precise statement of results}

Following \cite[Section 2]{GTV97} we first rewrite \eqref{CP} as the first order system
\begin{equation}
 \label{firstorder}
(i\partial_t \mp J_x)u_{\pm}=\mp \frac14 J_x^{-1}B_k(u_++u_-)\mp\frac12 J_x^{-1}(u_++u_-),
\end{equation}
where $J_x=(1-\Delta_x)^{\frac12}$ is the Bessel potential operator of order $-1$ in the space variable $x$, and
$u_{\pm}=u \pm iJ_x^{-1}\partial_tu$ so that the initial conditions become
\begin{equation}
 \label{data}
u_{\pm}(0)=u_0 \pm iJ_x^{-1}u_1=:f_{\pm}\in \widehat{H}^r_s(\R^3).
\end{equation}
To treat the system \eqref{firstorder} with data given by \eqref{data} we will use the function spaces $X^{r,\pm}_{s,b}$ defined by
their norm
$$\n{u}{X^{r,\pm}_{s,b}}:= \left(\int d \xi d \tau \langle \xi \rangle^{sr'}\langle \tau \pm |\xi|\rangle^{br'} |\F{u}(\xi , \tau)|^{r'} \right) ^{\frac{1}{r'}}.$$
For $s=b=0$ we write $\widehat{L^r_{xt}}:=X^{r,+}_{0,0}=X^{r,-}_{0,0}$, correspondingly we set $\widehat{L^r_{x}}:=\widehat{H}^r_0$.
Local solutions are obtained by the contraction mapping principle in the time restriction space
$$X^{r,\pm}_{s,b}(\delta) := \{u = \tilde{u}|_{[-\delta,\delta] \times \R^3} : \tilde{u} \in X^{r,\pm}_{s,b}\}$$
endowed with the norm
$$\n{u}{X^{r,\pm}_{s,b}(\delta)}:= \inf \{ \n{\tilde{u}}{X^{r,\pm}_{s,b}} : \tilde{u}|_{[-\delta,\delta] \times \R^3} =u\}  .$$
We will always have $b > \frac{1}{r}$, hence $X^{r,\pm}_{s,b}\subset C(\R,\widehat{H}^r_s )$ and 
$X^{r,\pm}_{s,b}(\delta)\subset C([-\delta,\delta],\widehat{H}^r_s )$, which gives the persistence property of our solutions. To deal with $B_2$ - especially if time derivatives are involved - we also need the norms
$$\|u\|_{Z^{r,\pm}_{s,b}}:=\|u\|_{X^{r,\pm}_{s,b}}+\|\partial_tu\|_{X^{r,\pm}_{s-1,b}};$$
the corresponding restriction spaces are defined precisely as above. Now our result concerning $B_1$ reads as follows.
\begin{theorem}
 \label{B1local}
Let $1<r\le 2$, $s>\frac{2}{r}$, $\frac{1}{r}<b<1$ and $f_{\pm} \in \widehat{H}^r_s$. Then there exist $\delta=\delta(\|f_+\|_{\widehat{H}^r_s},\|f_-\|_{\widehat{H}^r_s})>0$ and a unique solution 
$(u_+,u_-) \in X^{r,+}_{s,b}(\delta) \times X^{r,-}_{s,b}(\delta)$ of \eqref{firstorder} with $k=1$ satisfying
the initial condition \eqref{data}. The solution is persistent and the flow map
$$(f_+,f_-) \mapsto (u_+,u_-), \qquad \widehat{H}^r_s \times \widehat{H}^r_s \rightarrow 
X^{r,+}_{s,b}(\delta) \times X^{r,-}_{s,b}(\delta)$$
is locally Lipschitz continuous.
\end{theorem}
Similarly we can show the following for $B_2$.
\begin{theorem}
 \label{B2local}
Let $1<r\le 2$, $s>\frac{2}{r}+1$, $\frac{1}{r}<b<1$ and $f_{\pm} \in \widehat{H}^r_s$. Then there exist $\delta=\delta(\|f_+\|_{\widehat{H}^r_s},\|f_-\|_{\widehat{H}^r_s})>0$ and a unique solution 
$(u_+,u_-) \in Z^{r,+}_{s,b}(\delta) \times Z^{r,-}_{s,b}(\delta)$ of \eqref{firstorder} with $k=2$ satisfying
the initial condition \eqref{data}. The solution is persistent and the flow map
$$(f_+,f_-) \mapsto (u_+,u_-), \qquad \widehat{H}^r_s \times \widehat{H}^r_s \rightarrow 
Z^{r,+}_{s,b}(\delta) \times Z^{r,-}_{s,b}(\delta)$$
is locally Lipschitz continuous.
\end{theorem}
The general LWP theorem \cite[Theorem 2.3]{G04} reduces the proofs of Theorem \ref{B1local} and \ref{B2local} to that of bilinear estimates in $X^{r,\pm}_{s,b}$-norms. The next section is devoted to the proof of the following key estimate.
\begin{theorem}
 \label{key}
Let $1<r\le 2$, $b > \frac{1}{r}$, and $\sigma > \frac{2}{r}$. Then
\begin{equation}
 \label{keyest}
\|J_x^{\sigma}(uv)\|_{\widehat{L^r_{xt}}}+\|J_x^{\sigma-1}\partial_t(uv)\|_{\widehat{L^r_{xt}}}
\ls \|u\|_{X^{r,\pm}_{\sigma,b}}\|v\|_{X^{r,[\pm]}_{\sigma,b}},
\end{equation}
where $[\pm]$ denotes independent signs.
\end{theorem}
Assume \eqref{keyest} already to be proven. Concerning $B_1$ we then have that for all $b,r$ and $s=\sigma$ according to the
assumptions of Theorem \ref{key} and $b'\le 0$
\begin{multline*}
\|J_x^{-1}\partial(uv)\|_{X^{r,\pm}_{s,b'}}\le \|J_x^{s-1}\partial(uv)\|_{\widehat{L^r_{xt}}}\\
\le \|J_x^{s}(uv)\|_{\widehat{L^r_{xt}}}+\|J_x^{s-1}\partial_t(uv)\|_{\widehat{L^r_{xt}}}
\ls \|u\|_{X^{r,\pm}_{s,b}}\|v\|_{X^{r,[\pm]}_{s,b}},
\end{multline*}
which combined with \cite[Theorem 2.3]{G04} leads to Theorem \ref{B1local}, since the linear term on the right of \eqref{firstorder}
can be trivially taken care of. Similarly for $B_2$ we have with $s=\sigma + 1 > 1 + \frac{2}{r}$ and $r,b,b'$ as before
\begin{multline*}
\|J_x^{-1}(\partial u \partial v)\|_{Z^{r,\pm}_{s,b'}}
\le \|J_x^{\sigma}(\partial u \partial v)\|_{\widehat{L^r_{xt}}}+
\|J_x^{\sigma-1}\partial_t( \partial u \partial v)\|_{\widehat{L^r_{xt}}}\\
\ls \|\partial u\|_{X^{r,\pm}_{\sigma,b}}\|\partial v\|_{X^{r,[\pm]}_{\sigma,b}}
\ls \|u\|_{Z^{r,\pm}_{s,b}}\|v\|_{Z^{r,[\pm]}_{s,b}},
\end{multline*}
which is sufficient for Theorem \ref{B2local}.

\section{Proof of the key estimate}

Theorem \ref{key} will be a consequence of several bilinear estimates for free solutions of the half-wave equations
$(i\partial_t \pm D_x)u=0$, subject to the initial condition $u(0)=u_0$. So for the remaining part of the paper let
$u_{\pm}(t)=e^{\pm itD_x}u_0=\F_x^{-1}e^{\pm it |\xi|}\F_x u_0$ and $v_{\pm}(t)=e^{\pm itD_x}v_0$. By the tranfer principle
- see e. g. \cite[Proposition 3.5]{KS} or \cite[Lemma 2.1]{G04} - the proof of \eqref{keyest} essentially\footnote{for low
frequencies $|\xi|\le 1$ we obtain $\|J_x^{\sigma}(uv)\|_{\widehat{L^r_{xt}}}\ls \|u\|_{X^{r,\pm}_{\sigma,b}}\|v\|_{X^{r,[\pm]}_{\sigma,b}}$ by Young's inequality and Sobolev type embeddings.} reduces to showing
that
\begin{equation}
 \label{keyfree}
\|J_x^{\sigma -1}\partial_x (u_{\pm}v_{[\pm]})\|_{\widehat{L^r_{xt}}}+\|J_x^{\sigma -1}\partial_t (u_{\pm}v_{[\pm]})\|_{\widehat{L^r_{xt}}} \ls \|u_0\|_{\widehat{H}^r_{\sigma}}\|v_0\|_{\widehat{H}^r_{\sigma}}.
\end{equation}
To prove \eqref{keyfree} we make substantial use of the calculations in \cite{FK}. By symmetry it suffices to consider the $(++)$-
and $(+-)$-cases. For both we calculate the space-time Fourier transform of the product. Defining $P_{\pm}(\eta):=
|\frac{\xi}{2}-\eta|\pm|\frac{\xi}{2}+\eta|$ with $\nabla P_{\pm}(\eta)= \frac{\eta- \frac{\xi}{2}}{|\eta- \frac{\xi}{2}|}
\pm \frac{\eta+ \frac{\xi}{2}}{|\eta+ \frac{\xi}{2}|}$ and using the properties of the $\delta$-distribution we obtain
$$\F u_+v_{\pm}(\xi,\tau)=c \int_{P_\pm(\eta)=\tau} \frac{dS_{\eta}}{|\nabla P_{\pm}(\eta)|}
\widehat{u_0}(\frac{\xi}{2}-\eta)\widehat{v_0}(\frac{\xi}{2}+\eta),$$
for more details see \cite[Sections 3 and 4]{FK}. Observe that the set $\{P_+(\eta)=\tau \}$ ($\{P_-(\eta)=\tau \}$) is an
ellipsoid (hyperboloid) of rotation, so the $(++)$-case ($(+-)$-case) is henceforth referred to as elliptic (hyperbolic).

\subsection{The elliptic case}

We choose $0 < s_{1,2}< \frac{2}{r}$ with $s_1+s_2=\frac{2}{r}$ and use H\"older's inequality to get
\begin{multline*}
|\F u_+v_+(\xi,\tau)| \ls \left(\int_{P_+(\eta)=\tau} \frac{dS_{\eta}}{|\nabla P_{+}(\eta)|} 
|\tfrac{\xi}{2}-\eta|^{-s_1r}|\tfrac{\xi}{2}+\eta|^{-s_2r}\right)^{\frac{1}{r}} \times \\
\left(\int_{P_+(\eta)=\tau} \frac{dS_{\eta}}{|\nabla P_{+}(\eta)|}
|\widehat{J_x^{s_1}u_0}(\tfrac{\xi}{2}-\eta)\widehat{J_x^{s_2}v_0}(\tfrac{\xi}{2}+\eta)|^{r'}\right)^{\frac{1}{r'}}.
\end{multline*}
For the first factor we apply \cite[Lemma 4.1]{FK} to see that
\begin{multline*}
 \int_{P_+(\eta)=\tau} \frac{dS_{\eta}}{|\nabla P_{+}(\eta)|} 
|\tfrac{\xi}{2}-\eta|^{-s_1r}|\tfrac{\xi}{2}+\eta|^{-s_2r} \\
= c \int_{-1}^1 |\tau + |\xi|x|^{1-s_1r}|\tau - |\xi|x|^{1-s_2r}dx \\
= c \int_{-1}^1 |\tfrac{\tau}{|\xi|} + x|^{1-s_1r}|\tfrac{\tau}{|\xi|} - x|^{1-s_2r}dx \le c_{s_1,s_2}.
\end{multline*}
Taking the $L^{r'}_{\xi,\tau}$-norm of the second factor and using the coarea formula we arrive at
$$\|u_+v_+\|_{\widehat{L^r_{xt}}} \ls \|J^{s_1}_xu_0\|_{\widehat{L^r_{x}}} \|J^{s_2}_xv_0\|_{\widehat{L^r_{x}}}.$$
Unfortunately this argument breaks down, if $s_1=0$ or $s_2=0$, cf. the necessary condition (9) in \cite{FK}. To overcome this
difficulty we split $u_+v_+=P_{\gs}(u_+,v_+)+P_{\ll}(u_+,v_+)$, where
$$\F_xP_{\gs}(f,g)(\xi):=\int_{\xi_1+\xi_2=\xi}\widehat{f}(\xi_1)\widehat{g}(\xi_2)\chi_{\{|\xi_1|\gs|\xi_2|\}}d\xi_1.$$
By the preceeding we have
\begin{equation}
 \label{Helmut}
\|P_{\gs}(u_+,v_+)\|_{\widehat{L^r_{xt}}} \ls \|J^{\frac{2}{r}}_xu_0\|_{\widehat{L^r_{x}}} \|v_0\|_{\widehat{L^r_{x}}}.
\end{equation}
To estimate $P_{\ll}(u_+,v_+)$ we decompose $u_0$ dyadically into $u_0 = \sum_{k\ge 0}P_{\Delta k}u_0$ with
$P_{\Delta 0}= \F_x^{-1}\chi_{\{|\xi| \le 1\}}\F_x$ and, for $k \ge 1$, $P_{\Delta k}= \F_x^{-1}\chi_{\{|\xi| \sim 2^k\}}\F_x$, so that
$$\|P_{\ll}(u_+,v_+)\|_{\widehat{L^r_{xt}}} \le \sum_{k\ge 0}\|P_{\ll}(P_{\Delta k}u_+,v_+)\|_{\widehat{L^r_{xt}}}.$$
Now by \cite[Lemma 12.2]{FK} we have
\begin{equation}
 \label{Gränke}
\int_{P_\pm(\eta)=\tau} \frac{dS_{\eta}}{|\nabla P_{\pm}(\eta)|}\chi_{\{2^k \sim |\tfrac{\xi}{2}-\eta|\ll |\tfrac{\xi}{2}+\eta|\}}
\ls 2^{2k},
\end{equation}
hence a H\"older application as above gives
$$\|P_{\ll}(P_{\Delta k}u_+,v_+)\|_{\widehat{L^r_{xt}}} \ls 2^{\frac{2k}{r}}\|P_{\Delta k}u_0\|_{\widehat{L^r_{x}}} \|v_0\|_{\widehat{L^r_{x}}}.$$
Summing up the dyadic pieces and combining the result with \eqref{Helmut} we obtain for $\sigma > \frac{2}{r}$
\begin{equation}
 \label{Stamm}
\|u_+v_+\|_{\widehat{L^r_{xt}}} \ls \|J^{\sigma}_xu_0\|_{\widehat{L^r_{x}}} \|v_0\|_{\widehat{L^r_{x}}}.
\end{equation}
The convolution constraint $\xi = \xi_1 + \xi_2 = (\tfrac{\xi}{2}-\eta)+(\tfrac{\xi}{2}+\eta)$ implies
$\lb\xi\rb^{\sigma} \ls \lb \xi_1\rb^{\sigma} + \lb \xi_2\rb^{\sigma} = \lb\tfrac{\xi}{2}-\eta\rb^{\sigma}+ 
\lb\tfrac{\xi}{2}+\eta\rb^{\sigma}$, and we may exchange $u_0$ and $v_0$ in \eqref{Stamm}. This gives
\begin{equation}
 \label{Ober}
\|J_x^{\sigma}(u_+v_+)\|_{\widehat{L^r_{xt}}} \ls \|J_x^{\sigma}u_0\|_{\widehat{L^r_{x}}}\|J_x^{\sigma}v_0\|_{\widehat{L^r_{x}}},
\end{equation}
provided $\sigma > \frac{2}{r}$. In \eqref{Ober} we may clearly replace $J_x^{\sigma}(u_+v_+)$ by
$J_x^{\sigma - 1}\partial_t(u_+v_+)$, since in the support of $\F (u_+v_+)$ we have $\tau = |\tfrac{\xi}{2}-\eta|+ |\tfrac{\xi}{2}+\eta|$ and hence $\lb \xi \rb^{\sigma -1}|\tau| \le \lb\tfrac{\xi}{2}-\eta\rb^{\sigma}+ 
\lb\tfrac{\xi}{2}+\eta\rb^{\sigma}$. Thus we have shown:
\begin{lemma}
 \label{elliptic}
Let $1 \le r \le 2$ and $\sigma > \frac{2}{r}$. Then
$$\|J_x^{\sigma}(u_+v_+)\|_{\widehat{L^r_{xt}}}+\|J_x^{\sigma - 1}\partial_t(u_+v_+)\|_{\widehat{L^r_{xt}}}
\ls \|u_0\|_{\widehat{H}^r_{\sigma}}\|v_0\|_{\widehat{H}^r_{\sigma}}.$$
\end{lemma}

\subsection{The hyperbolic case}

The estimation in this case goes along similar lines as in section 3.1, as long as
\begin{equation}
 \label{fertig}
|\frac{\xi}{2}-\eta|+ |\frac{\xi}{2}+\eta| \le c_1 |\xi|.
\end{equation}
If \eqref{fertig} is fulfilled, we choose again $s_{1,2}\in (0,\frac{2}{r})$ with $s_1+s_2=\frac{2}{r}$ and obtain from
\cite[Lemma 4.4]{FK} that
\begin{multline*}
 \int_{P_-(\eta)=\tau} \frac{dS_{\eta}}{|\nabla P_{-}(\eta)|} 
|\tfrac{\xi}{2}-\eta|^{-s_1r}|\tfrac{\xi}{2}+\eta|^{-s_2r} \chi_{\{\eqref{fertig}\}}\\
= c \int_{1}^{c_1}|\tfrac{\tau}{|\xi|} + x|^{1-s_1r}|\tfrac{\tau}{|\xi|} - x|^{1-s_2r}dx \le c_{s_1,s_2},
\end{multline*}
which gives
$$\|u_+v_-\|_{\widehat{L^r_{xt}}} \ls \|J^{s_1}_xu_0\|_{\widehat{L^r_{x}}} \|J^{s_2}_xv_0\|_{\widehat{L^r_{x}}}$$
and hence
\begin{equation}
 \label{Ralf}
\|P_{\gs}(u_+,v_-)\|_{\widehat{L^r_{xt}}} \ls \|J^{\frac{2}{r}}_xu_0\|_{\widehat{L^r_{x}}} \|v_0\|_{\widehat{L^r_{x}}}.
\end{equation}
A dyadic decomposition together with \eqref{Gränke} shows that
\begin{equation}
 \label{Rolf}
\|P_{\ll}(P_{\Delta k}u_+,v_-)\|_{\widehat{L^r_{xt}}} \ls 2^{\frac{2k}{r}}\|P_{\Delta k}u_0\|_{\widehat{L^r_{x}}} \|v_0\|_{\widehat{L^r_{x}}},
\end{equation}
and combining \eqref{Ralf} and \eqref{Rolf} after summation in $k$ we arrive at
$$\|u_+v_-\|_{\widehat{L^r_{xt}}} \ls\|J^{\sigma}_xu_0\|_{\widehat{L^r_{x}}} \|v_0\|_{\widehat{L^r_{x}}},$$
provided $1 \le r \le 2$, $\sigma > \frac{2}{r}$, and $u_+v_-$ fulfills assumption \eqref{fertig}. To fix a partial result
concerning the hyperbolic case, let $P(u,v)$ denote the projection on the domain in Fourier space, where \eqref{fertig} holds.
Then, taking into account the arguments at the end of Section 3.1, we have the following estimate.
\begin{lemma}
 \label{hyperbolic1}
Let $1 \le r \le 2$ and $\sigma > \frac{2}{r}$. Then
$$\|J_x^{\sigma}P(u_+,v_-)\|_{\widehat{L^r_{xt}}}+\|J_x^{\sigma - 1}\partial_tP(u_+,v_-)\|_{\widehat{L^r_{xt}}}
\ls \|u_0\|_{\widehat{H}^r_{\sigma}}\|v_0\|_{\widehat{H}^r_{\sigma}}.$$
\end{lemma}
We turn to the region, where
\begin{equation}
 \label{fix}
c_1 |\xi| < |\frac{\xi}{2}-\eta|+ |\frac{\xi}{2}+\eta| ,
\end{equation}
and denote the projection thereon by $Q(u,v)$. We apply again \cite[Lemma 4.4]{FK} with $F(|\tfrac{\xi}{2}-\eta|,|\tfrac{\xi}{2}+\eta|)=|\tfrac{\xi}{2}-\eta|^{-s_1r}|\tfrac{\xi}{2}+\eta|^{-s_2r}\chi_{\{\eqref{fix}\}}$,
where $s_{1,2}\ge 0$ and $s_1+s_2 = \frac{3}{r}+\e$. This gives
\begin{multline*}
 \int_{P_-(\eta)=\tau} \frac{dS_{\eta}}{|\nabla P_{-}(\eta)|} 
|\tfrac{\xi}{2}-\eta|^{-s_1r}|\tfrac{\xi}{2}+\eta|^{-s_2r} \chi_{\{\eqref{fix}\}}\\
= c \int_{c_1}^{\infty}|\tau + |\xi|x|^{1-s_1r}|\tau - |\xi|x|^{1-s_2r}dx \\
= c |\xi|^{2-(s_1+s_2)r}\int_{c_1}^{\infty}|\tfrac{\tau}{|\xi|} + x|^{1-s_1r}|\tfrac{\tau}{|\xi|} - x|^{1-s_2r}dx \ls |\xi|^{2-(s_1+s_2)r},
\end{multline*}
which in turn implies
\begin{equation}
 \label{Heinrich}
\|D_x^{s_1+s_2-\frac{2}{r}}Q(u_+,v_-)\|_{\widehat{L^r_{xt}}} \ls
\|J_x^{s_1}u_0\|_{\widehat{L^r_{x}}}\|J_x^{s_2}v_0\|_{\widehat{L^r_{x}}}.
\end{equation}
Bilinear interpolation of \eqref{Heinrich} with $r=1$ and
$$\|u_+v_-\|_{L^2_{xt}}\ls \|J_x^{\sigma_1}u_0\|_{L_x^2}\|J_x^{\sigma_2}v_0\|_{L_x^2}, 
\qquad (\sigma_{1,2}\ge 0, \sigma_1+\sigma_2>1),$$
which follows from Strichartz estimate, gives the sharpened version
$$\|D_x^{s}Q(u_+,v_-)\|_{\widehat{L^r_{xt}}} \ls
\|J_x^{s_1}u_0\|_{\widehat{L^r_{x}}}\|J_x^{s_2}v_0\|_{\widehat{L^r_{x}}},$$
where $1 \le r \le 2$, $s=(1-\frac{2}{r'})(1+\e)$, $s_{1,2}\ge 0$ with $s_1+s_2= 3-\frac{4}{r'}+\e$ and $\e > 0$. If in addition $r>1$ and $\e$ is sufficiently small, so that $s \le 1$, we may replace the $D_x^s$ by $J_x^{s-1}\partial_x$ and hence by
$J_x^{-\frac{2}{r'}}\partial_x$. This gives
$$\|J_x^{-\frac{2}{r'}}\partial_xQ(u_+,v_-)\|_{\widehat{L^r_{xt}}} \ls
\|J_x^{s_1}u_0\|_{\widehat{L^r_{x}}}\|J_x^{s_2}v_0\|_{\widehat{L^r_{x}}}$$
for all $r \in (1,2]$ and $s_{1,2}\ge0$ with $s_1+s_2>3-\frac{4}{r'}$. Using once more $\lb\xi\rb \le 
\lb\tfrac{\xi}{2}-\eta\rb+\lb\tfrac{\xi}{2}+\eta\rb$ we conclude for $\sigma > \frac{2}{r}$ that
$$\|J_x^{\sigma - 1}\partial_xQ(u_+,v_-)\|_{\widehat{L^r_{xt}}}
\ls \|u_0\|_{\widehat{H}^r_{\sigma}}\|v_0\|_{\widehat{H}^r_{\sigma}},$$
which also holds true with $\partial_t$ instead of $\partial_x$, since we are in the hyperbolic case, where $|\tau|\le|\xi|$.
Summarizing we have:
\begin{lemma}
 \label{hyperbolic2}
Let $1 < r \le 2$ and $\sigma > \frac{2}{r}$. Then
$$\|J_x^{\sigma - 1}\partial_xQ(u_+,v_-)\|_{\widehat{L^r_{xt}}}+\|J_x^{\sigma - 1}\partial_tQ(u_+,v_-)\|_{\widehat{L^r_{xt}}}
\ls \|u_0\|_{\widehat{H}^r_{\sigma}}\|v_0\|_{\widehat{H}^r_{\sigma}}.$$
\end{lemma}Now the crucial estimate \eqref{keyfree} follows from the Lemmas \ref{elliptic}, \ref{hyperbolic1}, and \ref{hyperbolic2}.

\end{document}